\documentclass[12pt]{article}
\usepackage[utf8]{inputenc}
\usepackage[centertags]{amsmath}
\usepackage[hidelinks]{hyperref}
\usepackage{amsfonts}
\usepackage{amssymb}
\usepackage{latexsym}
\usepackage{amsthm}
\usepackage{newlfont}
\usepackage{graphicx}
\usepackage{listings}
\usepackage{booktabs}
\usepackage{abstract}
\usepackage{enumerate}
\usepackage{xcolor}
\usepackage[a4paper,margin=30mm]{geometry}

\RequirePackage{srcltx}
\date{}

\newcommand{\actaqed}{\hfill $\actabox$}
{\medskip\noindent \textit{Proof of #1. }}%
{\actaqed \medskip}

\def\cC{{\mathcal C}}

\def\bbR{{\mathbb R}}

\def\bbT{{\mathbb T}}

\def\bF{{\mathbf F}}

\def\bW{{\mathbf W}}

\def\bx{\mathbf x}
\def\by{\mathbf y}

 \def \<{\langle}
\def\>{\rangle}

\def \e{\varepsilon}

\def \ff{\varphi}

\def\fee{\varphi}

\newtheorem{Theorem}{Theorem}[section]
\newtheorem{Lemma}{Lemma}[section]

\newtheorem{Remark}{Remark}[section]

\newtheorem{Corollary}{Corollary}[section]
\numberwithin{equation}{section}

\newtheorem{exttheorem}{Theorem}[section]

\newtheorem{extlemma}{Lemma}[section]

\newcommand{\be}{\begin{equation}}
\newcommand{\ee}{\end{equation}}

\begin{document}

\title{A remark on the entropy numbers of mixed smoothness classes}

\author{Yu. Malykhin and V. Temlyakov}

\newcommand{\Addresses}{{
  \bigskip
  \footnotesize
  
Yu.V.~Malykhin, \textsc{Steklov Mathematicl Institute of Russian Academy of
Sciences, Russia, Lomonosov Moscow State University, Russia.\\
E-main:} \texttt{malykhin-yuri@yandex.ru}

\medskip
V.N. Temlyakov, \textsc{University of South Carolina, USA,\\ Steklov Mathematical Institute of Russian Academy of Sciences, Russia,\\ Lomonosov Moscow State University, Russia, \\ Moscow Center of Fundamental and Applied Mathematics, Russia.\\  
E-mail:} \texttt{temlyakovv@gmail.com}

}}
\maketitle

\begin{abstract}{This paper is a followup paper to  the very recent breakthrough paper of Bilyk, Neumayer, and Ullrich. We focus on 
the entropy numbers of classes of multivariate functions with bounded mixed derivative. Basically, this paper completes the study 
of the asymptotic behavior of the entropy numbers of the corresponding classes in the case of large smoothness. 
In particular, we find the right orders of both the entropy numbers and the Kolmogorov widths of classes of multivariate functions with bounded mixed derivative, when the class is defined by the restriction in the uniform norm and the error of approximation is measured in the uniform norm. The case of functions of two variables was settled about 30 years ago and the case of three or more variables stayed open. } 
\end{abstract}

\section{Introduction}
\label{In}

In this paper we study asymptotic behavior of the entropy numbers of function classes defined by restrictions imposed on the corresponding mixed derivative of a function. These are classes from the family of function classes with mixed smoothness. Systematic study of asymptotic characteristics of these classes began at the end of 1950th. It was understood in the very beginning that such a study is a difficult task. Some problems are still open, especially in the case of small smoothness. It is well known that the study of numerical integration of these classes is closely connected with the discrepancy problems (see, for instance, \cite{VT170}) and the study of the entropy numbers of these classes is closely connected with 
an outstanding problem from the probability theory, namely, the Small Ball Problem (see, for instance, \cite{KL} and \cite{VTbook}, Ch.3). The reader can find results on approximation of functions with mixed smoothness and related results in the books \cite{VTmon}, \cite{DTU}, and \cite{VTbookMA}. 

This paper is motivated by the very recent 
paper \cite{BNU}, where the Small Ball Problem was solved via studying the entropy numbers of classes with mixed smoothness. 
First, we briefly formulate the known results and main results of the paper, then we give  detailed historical comments in Section \ref{H}.
Section~\ref{sec_main} contains proofs of new results. Section \ref{D} is devoted to a discussion. 

Let $X$ be a Banach space and let $B_X$ denote the unit ball of $X$ with the center at $0$. Denote by $B_X(y,r)$ a ball with center $y$ and radius $r$: $\{x\in X:\|x-y\|\le r\}$. For a compact set $A$ and a positive number $\e$ we define the covering number $N_\e(A)$:
$$
N_\e(A) := N_\e(A,X) :=\min \{n : \exists y^1,\dots,y^n :A\subseteq \cup_{j=1}^n B_X(y^j,\e)\}
$$
and the entropy numbers $\e_k(A,X)$:
$$
\e_k(A,X) := \inf \{\e : \exists y^1,\dots ,y^{2^k} \in X : A \subseteq \cup_{j=1}
^{2^k} B_X(y^j,\e)\}.
$$
We refer the reader to \cite{VTbookMA} for detailed information on the entropy numbers of mixed smoothness classes. 
In this paper we discuss the case, when $X$ is the $L_p$ space, $1\le p\le \infty$. Let $(\Omega,\mu)$ be a probability space. 
  By the $L_p$, $1\le p< \infty$, norm we understand
$$
\|f\|_p:=\|f\|_{L_p(\Omega,\mu)} := \left(\int_\Omega |f|^p\,d\mu\right)^{1/p}.
$$
By the $L_\infty$-norm we understand the uniform norm of continuous functions
$$
\|f\|_\infty := \sup_{\omega\in\Omega} |f(\omega)|
$$
and with some abuse of notation we write $L_\infty(\Omega)$ for the space $\cC(\Omega)$ of continuous functions on $\Omega$. 
We consider multivariate $2\pi$-periodic (in each variable) functions defined on $\Omega = \bbT^d$ equipped with the normalized 
Lebesgue measure $d\mu = (2\pi)^{-d}d\bx$. 

Let $F_r(x,\alpha)$ be the univariate Bernoulli kernels 
$$
F_r(x,\alpha) := 1+2\sum_{k=1}^\infty k^{-r}\cos (kx-\alpha \pi/2).
$$
 For $\bx=(x_1,\dots,x_d)$ and $\alpha=(\alpha_1,\dots,\alpha_d)$ we define
$$
F_r({\bx},\alpha):= \prod_{i=1}^d F_r(x_i,\alpha_i)
$$
and
$$
\bW^r_{q,\alpha}:=\{f:f=F_r(\cdot,\alpha)\ast \fee,\quad \|\fee\|_q\le 1\},
$$
where
$$
( F_r(\cdot,\alpha) \ast \varphi)(\bx):= (2\pi)^{-d}\int_{\bbT^d} F_r(\bx-\by,\alpha) \ff(\by)d\by.
$$
where $\ast$ means convolution. In the univariate case we use the notation $W^r_{q,\alpha}$.

The classes $\bW^r_{q,\alpha}$ are classical classes of functions with {\it dominating mixed derivative}
(Sobolev-type classes of functions with mixed smoothness).

In this paper we study the asymptotic behavior of $\e_k(\bW^r_q,L_p)$. We use 
the standard notation for the order inequalities. For two sequences $\{a_k\}_{k=1}^\infty$ and $\{b_k\}_{k=1}^\infty$ we write $a_k\ll b_k$, which is the same as $b_k\gg a_k$, if there exists a constant $c>0$ such that $a_k\le cb_k$ for all $k\ge 1$. We write $a_k\asymp b_k$ if $a_k\ll b_k$ and $b_k\ll a_k$. 

We now formulate our new results and very briefly mention some known results. The reader can find the detailed historical comments in Section \ref{H}. We only discuss here the case of large smoothness, when $r > r(q) := \max(1/q,1/2)$. For convenience we 
use the notation
$$
k(d) := k (\log(k+1))^{1-d},
\quad k=1,2,\ldots.
$$
We begin with results, which were known for long time (for details see 
Section \ref{H} below) for all $d\ge 2$. The following relation 
\begin{equation}\label{InH2}
\e_k(\bW^r_{q,\alpha},L_p)\asymp k(d)^{-r}  
\end{equation}
holds for pairs $(q,p)$ coming from the following sets
\be\label{Sets1}
S_1:= \{(q,p): 1<q,p<\infty\},\quad S_2:= \{(q,p): 1<q\le \infty, p=1\},
\ee
\be\label{Sets1a}
S_3:= \{(q,p): q=\infty, 1\le p<\infty\}.
\ee
Next, we formulate new results, in particular, very recent results from the breakthrough paper \cite{BNU}, which covers 
the case $1<q<\infty$, $p=\infty$. Note that all these results were known in the case $d=2$ (for details see 
Section \ref{H} below). The following relation
\begin{equation}\label{InH2a}
\e_k(\bW^r_{q,\alpha},L_p)\asymp k(d)^{-r}  (\log(k+1))^{1/2}
\end{equation}
holds for pairs $(q,p)$ coming from the following sets
\be\label{Sets2}
  S_4:= \{(q,p): q=1, 1\le p<  \infty \},\quad S_5:= \{(q,p): 1< q\le \infty, p=\infty\}.
\ee
In the case $(q,p) \in S_6:=\{(1,\infty)\}$ we present the following result
\begin{equation}\label{InH2b}
\e_k(\bW^r_{1,0},L_\infty)\asymp k(d)^{-r}  \log(k+1).
\end{equation}

Note that the case $(q,p) = (\infty,\infty)$ is covered in the previous version  \cite{MT}, v1 of this paper.

 \section{Historical comments}
\label{H}

It is well known that in the univariate case  (see \cite{BS})
\begin{equation}\label{H1}
\e_k(W^r_{q,\alpha},L_p)\asymp k^{-r}
\end{equation}
 holds for all $1\le q,p \le \infty$ and $r>(1/q-1/p)_+$. We note that condition $r>(1/q-1/p)_+$ is a necessary and sufficient condition for compact embedding of $W^r_{q,\alpha}$ into $L_p$. Thus, (\ref{H1}) provides a complete description of the rate of $\e_k(W^r_{q,\alpha},L_p)$ in the univariate case. We point out that (\ref{H1}) shows that the rate of decay of $\e_k(W^r_{q,\alpha},L_p)$ depends only on $r$ and does not depend on $q$ and $p$. In this sense the strongest upper bound (for $r>1$) is $\e_k(W^r_{1,\alpha},L_\infty) \ll k^{-r}$ and the strongest lower bound is $\e_k(W^r_{\infty,\alpha},L_1)\gg k^{-r}$. 

There are different generalizations of classes $W^r_{q,\alpha}$ to the case of multivariate functions. In this section we only discuss known results for classes $\bW^r_{q,\alpha}$ of functions with bounded mixed derivative. For further discussions see \cite{VTbook}, Chapter 3 and \cite{DTU}.    

The following two theorems are from \cite{TE1} and \cite{TE2}.

\begin{exttheorem}\label{Wup} For $r>1$ and $1<q,p<\infty$ one has
$$
\e_k(\bW^r_{q,\alpha},L_p) \ll k(d)^{-r}.
$$
\end{exttheorem}

\begin{exttheorem}\label{Wlo} For $r>0$ and $1\le q <\infty$ one has
$$
\e_k(\bW^r_{q,\alpha},L_1) \gg k(d)^{-r}.
$$
\end{exttheorem}

The problem of estimating $\e_k(\bW^r_{q,\alpha},L_p)$ has a long history. The first result on the right order of $\e_k(\bW^r_{2,\alpha},L_2)$ was obtained by Smolyak \cite{Smo}. Later (see \cite{TE1}, \cite{TE2} and theorems above) it was established that
\begin{equation}\label{H2}
\e_k(\bW^r_{q,\alpha},L_p)\asymp k(d)^{-r}  
\end{equation}
holds for all $1<q,p<\infty$, $r>1$. The case $1<q=p<\infty$, $r>0$ was
established by Dinh Dung \cite{DD}. Belinsky \cite{Bel} extended (\ref{H2}) to the case $r>(1/q-1/p)_+$ when $1<q,p<\infty$. 

It is known in approximation theory (see \cite{VTbookMA}) that investigation of asymptotic characteristics of classes $\bW^r_{q,\alpha}$ in $L_p$ becomes more difficult when $q$ or $p$ takes value $1$ or $\infty$ than when $1<q,p<\infty$. It turns out to be the case for $\e_k(\bW^r_{q,\alpha},L_p)$ too. It was discovered that in some of these extreme cases ($q$ or $p$ equals $1$ or $\infty$) relation (\ref{H2}) holds and in other cases it does not hold. We describe the picture in detail. It was proved in \cite{TE2} that (\ref{H2}) holds for $p=1$, $1<q<\infty$, $r>0$. It was also proved that (\ref{H2}) holds for $p=1$, $q=\infty$ (see \cite{Bel} for $r>1/2$ and  \cite{KTE2} for $r>0$). Summarizing, we state that (\ref{H2}) holds for $1<q,p<\infty$ and $p=1$, $1<q\le\infty$ for all $d$ (with appropriate  restrictions on $r$). This easily implies that (\ref{H2}) also holds for $q=\infty$, $1\le p<\infty$. 

In the case $d=2$ this problem was essentially solved more than 20 years ago. We now cite the corresponding results in the case $d=2$. The first result on the right order of $\e_k(\bW^r_{q,\alpha},L_p)$ in the case $p=\infty$ was obtained by Kuelbs and Li \cite{KL} for $q=2$, $r=1$. It was proved in \cite{TE3} that
\begin{equation}\label{H3}
\e_k(\bW^r_{q,\alpha},L_\infty)\asymp k^{-r}(\log k)^{r+1/2}
\end{equation}
holds for $1<q<\infty$, $r>1$. We note that the upper bound in (\ref{H3}) was
proved under condition $r>1$ and the lower bound in (\ref{H3}) was proved under
condition $r>1/q$. Belinsky \cite{Bel} proved the upper bound in (\ref{H3}) for $1<q<\infty$ under condition $r>\max(1/q,1/2)$.   Relation (\ref{H3}) for $q=\infty$ under assumption $r>1/2$ was proved in  \cite{TE4}. 

The case $q=1$, $1\le p\le \infty$ was settled by Kashin and Temlyakov \cite{KaTe03}. The authors proved that 
\begin{equation}\label{H3'} 
\e_k(\bW^r_{1,\alpha},L_p)\asymp k^{-r}(\log k)^{r+1/2}
\end{equation}
holds for $1\le p<\infty$, $r>\max(1/2,1-1/p)$ and
\begin{equation}\label{H4} 
\e_k(\bW^r_{1,0},L_\infty)\asymp k^{-r}(\log k)^{r+1},\quad r>1.
\end{equation}

Let us make an observation on the base of the above discussion. In the univariate case the entropy numbers $\e_k(W^r_{q,\alpha},L_p)$ have the same order of decay with respect to $k$ for all pairs $(q,p)$, $1\le q,p\le\infty$. In the case $d=2$ we have three different orders of decay of $\e_k(\bW^r_{q,\alpha},L_p)$ which depend on the pair $(q,p)$. For instance, in the case $1<q,p<\infty$ it is $k^{-r}(\log k)^r$, in the case $q=1$, $1<p<\infty$, it is $k^{-r}(\log k)^{r+1/2}$ and in the case $q=1$, $p=\infty$ it is $k^{-r}(\log k)^{r+1}$.

An analogue of the upper bound in (\ref{H3}) for any $d$ was obtained by
Belinsky \cite{Bel}: for $q>1$ and $r>\max(1/q,1/2)$ we have
\begin{equation}\label{H5} 
\e_k(\bW^r_{q,\alpha},L_\infty)\ll k(d)^{-r}(\log (k+1))^{1/2}.
\end{equation}

  Kuelbs and Li \cite{KL} discovered the fact that there is a tight relationship between Small Ball Problem  and the behavior of the entropy $\log N_\e(\bW^1_{2,\alpha},L_\infty) $. Based on results obtained  by Livshits and Tsirelson \cite{LiTs}, by Bass \cite{Bass}, and  by Talagrand \cite{Tal} for the Small Ball Problem, they proved in the case of $d=2$
\begin{equation}\label{H6}
\e_k(\bW^1_{2,\alpha},L_\infty) \asymp k^{-1}(\log k)^{3/2}.
\end{equation}
Proof of the most difficult part of (\ref{H6}) -- the lower bound -- is based on a special inequality, known now as the Small Ball Inequality, for the Haar polynomials proved by Talagrand \cite{Tal} (see \cite{TE5} for a simple proof). 

Very recently, Bilyk, Neumayer, and Ullrich \cite{BNU} have made an outstanding progress by proving the following lower bound. Let 
$d\ge 3$, $r>1/2$, and $2\le q<\infty$. Then
\be\label{H7}
\e_k(\bW^r_{q,\alpha},L_\infty) \gg k(d)^{-r}(\log (k+1))^{1/2}.
\ee
This lower bound combined with (\ref{H5}) gives the right order under listed above restrictions
\be\label{H7a}
\e_k(\bW^r_{q,\alpha},L_\infty) \asymp k(d)^{-r}(\log (k+1))^{1/2}.
\ee
This bound allowed the authors of \cite{BNU} to solve the Small Ball Problem for $d\ge 3$. They point out that the results of \cite{BNU} were obtained with the assistance of the AI. 


\section{Bounds on the entropy numbers}
\label{sec_main}

\subsection{Preliminary results}

It will be convenient for us to formulate some inequalities
in terms of the entropy numbers of operators. Let $S$ be a linear operator from
$X$ to $Y$. We define the $n$th entropy number of $S$ as $$
\e_n(S:X\to Y) := \e_n(S(B_X),Y)
$$
where $S(B_X)$ is the image of the unit ball $B_X$ of $X$ under mapping $S$.

The entropy numbers of mixed smoothness class correspond to the entropy
numbers of the operator 
$$
I^r_\alpha \ff := F_r(\cdot,\alpha) \ast \ff(\cdot).
$$
Namely, we have
$$
\e_k(\bW^r_{q,\alpha},L_p) = \e_k(I^r_\alpha\colon L_q\to L_p).
$$

\paragraph{Interpolation.}
The following interpolation-type results are known (see also \cite{VT155}).

\begin{exttheorem}[\cite{VTbookMA}, Theorem 7.3.3]
\label{PT1}
For any $1\le p<v<q \le \infty$ and any Banach space $Y$ we have
$$
\e_{n+m}(S\colon L_v\to Y)\le 2\,\e_n(S\colon L_p\to Y)^{1-\theta}\e_m(S\colon L_q\to Y)^{\theta},
$$
where $1/v=(1-\theta)/p + \theta/q$.
\end{exttheorem}

\begin{extlemma}[\cite{VTbookMA}, Remark 7.3.2]
\label{OL4}
Let $1\le p<v<q\le\infty$ and let $A\subset L_p\cap L_q$. Then $A\subset L_v$ and
\be\label{O7}
\e_{n+m}(A,L_v) \le 2\,\e_n(A,L_p)^{1-\theta}\e_m(A,L_q)^{\theta},
\ee
where $1/v=(1-\theta)/p + \theta/q$.
\end{extlemma}

\paragraph{Duality.}
For a Banach space $X$ denote $X^*$ its dual and for $F\in X^*$ and $f\in X$ we  write for convenience
 $$
 \<f,F\> := F(f).
 $$
 For two Banach spaces $X$, $Y$ and a bounded linear operator $A:\, X\to Y$ the dual operator 
$A^*:\, Y^*\to X^*$ is defined by  
$$
\<Ax,y^*\> = \<x,A^*y^*\>
$$
for all   $x\in X$ and all $y^* \in Y^*$.

The following entropy duality result is a corollary of Theorem 5 from
\cite{AMST} (see also~\cite{Mil07}).

\begin{exttheorem}
Let $1<p<\infty$ and let $S$ be a compact operator from a Banach space $X$ to an
$L_p$-space.
Then there exist two constants $c=c(p)$ and $C=C(p)$ such that 
    $$
    C^{-1}\e_{ck}(S^*\colon L_{p'} \to X^*)
    \le \e_k(S\colon X \to L_p)
    \le C \e_{k/c}(S^*\colon L_{p'}\to X^*).
    $$
\end{exttheorem}

\begin{Corollary}
\label{OT2}
    Let $1<p<\infty$, $1\le q\le\infty$, and $r>(1/q-1/p)_+$. Then

\be\label{O5}
    \e_{ck}(\bW^r_{p',-\alpha}, L_{q'}) \le C \e_k(\bW^r_{q,\alpha}, L_p),
\ee
\be\label{O6}
    \e_{ck}(\bW^r_{q,\alpha}, L_p)
    \le C\e_{k}(\bW^r_{p',-\alpha}, L_{q'}).
\ee
\end{Corollary}

\begin{proof}
    Indeed, take $S:=I^r_\alpha$ and use that $S^*=I^r_{-\alpha}$.
\end{proof}

\subsection{The case $q=p=\infty$}
 
The main result of this subsection is the following new lower bound.

\begin{Theorem}\label{InT1} Let $d\ge 3$ and $r>1/2$. Then
\be\label{In1}
    \e_k(\bW^r_{\infty,\alpha},L_\infty) \gg k(d)^{-r}(\log (k+1))^{1/2}.
\ee
\end{Theorem}

The lower bound (\ref{In1}) combined with the known upper bound~\eqref{H5} gives the right order
\be\label{In1a}
\e_k(\bW^r_{\infty,\alpha},L_\infty) \asymp k(d)^{-r}(\log (k+1))^{1/2}, \quad r>1/2.
\ee

\begin{proof}[Proof of Theorem~\ref{InT1}]
The proof is based on the relation
(\ref{H7a}) and on the trick, which was used in the paper \cite{TE4}.
Let us use Theorem \ref{PT1} with operator $S = I^r_\alpha$ and
spaces $L_2$, $L_4$, $L_\infty$:
\be\label{P4}
\e_{2k}(\bW^r_{4,\alpha},L_\infty)\le 2\e_k(\bW^r_{2,\alpha},L_\infty)^{1/2} \e_k(\bW^r_{\infty,\alpha},L_\infty)^{1/2}.
\ee
From the the lower bound (\ref{H7}) we obtain
\be\label{P5}
\e_{2k}(\bW^r_{4,\alpha},L_\infty) \gg k(d)^{-r}(\log (k+1))^{1/2}.
\ee
The upper bound in (\ref{H5})  gives for $r>1/2$
\be\label{P6}
\e_{k}(\bW^r_{2,\alpha},L_\infty) \ll k(d)^{-r}(\log (k+1))^{1/2}.
\ee
Substituting bounds (\ref{P5}) and (\ref{P6}) into (\ref{P4}) we get the required lower bound.
\end{proof}

\subsection{The case $q=1$ and $1\le p<\infty$}
\label{subsec_q1} 

The main result of this subsection is the following Theorem.

\begin{Theorem}
\label{thm_q1}
Let $d\ge 2$, $1\le p<\infty$ and $r>\max(1/2,1-1/p)$. Then 
\be
    \label{O1}
    \e_k(\bW^r_{1,\alpha},L_p) \asymp k(d)^{-r} (\log(k+1))^{1/2}.
\ee
\end{Theorem}

\begin{proof}[Proof of Theorem~\ref{thm_q1}]
    For brevity we will use the notation $W:= \bW^r_{1,\alpha}$.

    First we prove the upper bound. Note that the right hand side of~\eqref{O1}
    does not depend on $p$; only the restriction on $r$ depends on it:
    $r>\max(1/2,1-1/p)$.
    Therefore, it is sufficient to 
    consider the case $2\le p<\infty$. By (\ref{O6}) and by (\ref{H5}) with
    $q=p'$ we obtain the required inequality
    \be\label{O2}
    \e_k(W,L_p)
    \ll \e_{c_1k}(\bW^r_{p',-\alpha},L_\infty)
    \ll k(d)^{-r} (\log(k+1))^{1/2}.
    \ee

    For the lower bound we first consider the case of the $L_2$ metric.
    Apply~\eqref{O5} with $p=2$ and the lower bound (\ref{H7}) with $q=2$ and
    obtain
    \be\label{O4}
    \e_k(W,L_2) \gg \e_{c_2k}(\bW^r_{2,-\alpha},L_\infty)
    \gg k(d)^{-r} (\log(k+1))^{1/2},
    \quad r>1/2.
    \ee
    Next, we use Lemma~\ref{OL4} with $m=n=k$, $A=W$, spaces $L_1$, $L_2$ and
    $L_q$, where $q>2$ is sufficiently close to $2$ (so that $r>1-1/q$):
    $$
    \e_{2k}(W,L_2) \le 2\e_k(W,L_1)^{1-\theta}\e_k(W,L_q)^{\theta}.
    $$
    The lower bound for the $\e_{2k}(W,L_2)$ and the upper bound for
    $\e_k(W,L_q)$ are already proved in~\eqref{O4} and~\eqref{O2}.
    Therefore, we obtain the required lower bound in $L_1$ (and for all
    $p$ by monotonicity).
\end{proof}

Theorem~\ref{thm_q1} covers the whole range $1\le p<\infty$ but the case $p=\infty$ is excluded. In the case $p=\infty$, $d\ge 2$ alike 
the case $p=\infty$, $d = 2$ (see~\eqref{H4}) we can treat the case $\alpha =0$. Here is the corresponding result.

\begin{Theorem}\label{OT3} Let $d\ge 2$  and $r>1$. Then 
\be\label{O8}
    \e_k(\bW^r_{1,0},L_\infty) \asymp k(d)^{-r} (\log(k+1)) .
\ee
\end{Theorem}
\begin{proof} We begin with the upper bounds. We use the following simple inequality, which easily follows from the definition of entropy numbers. For compact operators $A: X\to Y$ and $B: Y\to Z$ acting in Banach spaces $X$, $Y$, $Z$ we have
\be\label{O9}
\e_{n+m}(AB: X\to Z) \le \e_n(A: X\to Y)\e_m(B: Y\to Z).
\ee
For $r>1$ set $s:= r/2 >1/2$ and consider the operators $B:= I^s_0: L_1\to L_2$, $A:= I^s_0: L_2\to L_\infty$. Then 
$AB = I^r_0 : L_1\to L_\infty$. We now apply Theorem \ref{thm_q1} for
    the operator $A$, the bound (\ref{H5}) for the operator $B$
    and~\eqref{O9}:
    \begin{multline*}
    \e_{2k}(\bW^r_{1,0},L_\infty)
    = \e_{2k}(I^r_0\colon L_1\to L_\infty)
    \le \e_k(I^s_0\colon L_1\to L_2)\,\e_k(I^s_0\colon L_2\to L_\infty) = \\
    = \e_k(\bW^{s}_{1,0},L_2)\,\e_k(\bW^s_{2,0},L_\infty)
        \ll k(d)^{-s} (\log(k+1))^{1/2} \cdot k(d)^{-s}(\log(k+1))^{1/2}.
    \end{multline*}
    This implies the upper bound in~\eqref{O8}.

\begin{Remark}\label{OR1} The above proof works in the case of an arbitrary $\alpha$ as well. 
Indeed, we can take $B$ as above and $A:= I^s_\alpha: L_2\to L_\infty$.
\end{Remark}

We now proceed to the lower bounds. For that we need the following simple general Lemma \ref{OL5}. 

\begin{Lemma}\label{OL5} Let $E$ be a Banach space and $H$ be a Hilbert space and $A: E\to H$ be a compact operator. Then
\be\label{O10}
\e_k(A: E\to H)^2 \le 4 \e_k(A^*A: E\to E^*).
\ee
\end{Lemma}
\begin{proof} Denote $\e:= \e_k(A^*A: E\to E^*)$. By the definition, there exist
    $2^k$ points $v_j \in E^*$ such that for any 
$u\in B(E)$ there is $j\in \{1,2,\dots,2^k\}$ such that $\|A^*Au-v_j\|_{E^*} \le \e$.
    For each $j$ we can replace $v_j$ by $A^*Au_j$ with some $u_j\in B(E)$
    providing slighlty worse inequality: $\|A^*Au-A^*Au_j\|_{E^*}\le 2\e$.
    Therefore
$$
\|A(u-u_j)\|_H^2 = \<A^*A(u-u_j),u-u_j\> \le 4\e,
$$
    which means that the set $\{Au_j\}$ provides the required covering for the image
    $AB(E)$. This completes the proof of lemma.
\end{proof}

We take $E=L_1$, $H=L_2$ and $A:= I^s_0: L_1\to L_2$ with $s=r/2$. Then $A^*A =
    I^r_0$ and we obtain
$$
    \e_k(\bW^r_{1,0}, L_\infty)
    = \e_k(I^r_0\colon L_1\to L_\infty)
    \ge \frac14\e_k^2(I^s_0\colon L_1\to L_2)
    = \frac14\e_k^2(\bW^s_{1,0}, L_2).
$$
It remains to use Theorem~\ref{thm_q1}.
\end{proof}

\section{Discussion}
\label{D}

We begin with presenting some results on the Kolmogorov widths of smoothness classes $\bW^r_{\infty,\alpha}$. Let $X$ be a Banach space and $\bF\subset X$ be a  compact subset of $X$. The quantities  (see \cite{Ko36})
$$
d_n (\bF, X) :=  \inf_{\{u_i\}_{i=1}^n\subset X}
\sup_{f\in \bF}
\inf_{c_i} \left \| f - \sum_{i=1}^{n}
c_i u_i \right\|_X, \quad n = 1, 2, \dots,
$$
are called the {\it Kolmogorov widths} of $\bF$ in $X$.

We now illustrate a powerful method of proving lower bounds for the Kolmogorov widths. 
This method is based on the lower bounds for the entropy numbers.  For further discussion of related results on the entropy numbers we refer the reader to \cite{VTbook}, Ch.3, \cite{DTU}, and \cite{VTbookMA}, Ch.7. 
  The following lemma (see \cite{TE4} and \cite{VTbookMA}, p.229) plays a fundamental role in that method. Note that \cite{C} also can be used. 
\begin{Lemma}\label{DL1} Let $A$ be centrally symmetric compact set in a Banach space $X$ and for two real numbers $r>0$ and $a \in \bbR$ we have
$$
d_n(A,X) \ll n^{-r}(\log n)^a
$$
and 
$$
\e_k(A,X) \gg  k^{-r}(\log k)^a.
$$
Then the following relations
$$
d_n(A,X) \asymp \e_n(A,X) \asymp  n^{-r}(\log m)^a
$$ 
hold.
\end{Lemma}

For convenience we use the notation
$$
n(d) := n (\log (n+1))^{1-d}. 
$$
E. Belinsky (see  \cite{TrBe}) proved the following bounds
 \be\label{KolWb}
d_n(\bW^r_{2},L_\infty) \ll n(d)^{-r}(\log (n+1))^{1/2}, \quad r>1/2.
\ee

Theorem \ref{InT1}, Lemma \ref{DL1}, and the bound (\ref{KolWb}) imply the following result.

\begin{Theorem}\label{DT1} Let $d\ge 3$ and $r>1/2$. Then
\be\label{D2}
d_n(\bW^r_{\infty,\alpha},L_\infty) \asymp n(d)^{-r}(\log (n+1))^{1/2}.
\ee
\end{Theorem}

Note that an analog of Theorem \ref{DT1} for $d=2$ was proved in \cite{TE4}. 

We do not present a detailed historical discussion of the development of results on the Kolmogorov widths $d_n(\bW^r_{q,\alpha},L_p)$, $1\le q,p\le \infty$. The reader can find the corresponding comments in the books \cite{DTU} and \cite{VTbookMA}. For the very recent results we refer the reader to the paper \cite{BNU}. 

{\bf Conclusion.} Development of the theory of approximation of functions with mixed smoothness (other name for this theory is 
{\it hyperbolic cross approximations}) remains very important and some really difficult problems have been solved recently. 
The main contributions of this paper are the lower bounds in the case $d\ge 3$
\be\label{D3}
\e_k(\bW^r_{\infty,\alpha},L_\infty) \gg k(d)^{-r}(\log (k+1))^{1/2},
\ee
\be\label{D3a}
 d_n(\bW^r_{\infty,\alpha},L_\infty) \gg n(d)^{-r}(\log (n+1))^{1/2}.
\ee
These lower bounds complement and strengthen the following lower bounds from \cite{BNU} ($d\ge 3$, $q<\infty$)
\be\label{D4}
\e_k(\bW^r_{q,\alpha},L_\infty) \gg k(d)^{-r}(\log (k+1))^{1/2},
\ee
\be\label{D4a}
 d_n(\bW^r_{q,\alpha},L_\infty) \gg n(d)^{-r}(\log (n+1))^{1/2}.
\ee
We point out that our new bounds (\ref{D3}) and (\ref{D3a}) are based on the lower bounds (\ref{D4}) and (\ref{D4a}).

Theorem \ref{DT1} together with results from \cite{BNU} solve Open Problem 5.2
from \cite{VTbookMA}, p.243 under assumption $r>1/2$. Relation (\ref{In1a})
together with results from \cite{BNU} solve Open Problem 7.2 from
\cite{VTbookMA}, p.380 under assumptions $r>\max\{1/q,1/2\}$ and $1<q\le\infty$.
Theorem~\ref{thm_q1} solves Open Problem 7.1 from \cite{VTbookMA}, p.380 under assumptions $r>\max(1/2,1-1/p)$ and $1\le p<\infty$.

{\bf Acknowledgments.} Results of the subsection~\ref{subsec_q1} were obtained
with the assistance of ChatGPT, Astra. The authors have verified those results
and take full responsibility for their correctness.

 \Addresses
 
\end{document}